\documentclass{article}

\usepackage{amssymb, amsmath}

  \newtheorem{lemma}{Lemma}%

\newtheorem{remark}{Remark}

\newcommand{\vep}{\varepsilon}

\newcommand{\QED}{\hfill
$\underline{\underline{QED}}$}

\newcommand{\la}{\lambda}

\newcommand{\beq}{\begin{equation}}
\newcommand{\eeq}{\end{equation}}
\newcommand{\beqn}{\begin{eqnarray}}
\newcommand{\eeqn}{\end{eqnarray}}
\newcommand{\beqnn}{\begin{eqnarray*}}
\newcommand{\eeqnn}{\end{eqnarray*}}
\newcommand{\nn}{|\!|}

\newcommand{\rr}{\mathbb{R}}
\newcommand{\n}{\mathbb{N}}

\begin{document}
\title{On Neumann ``superlinear'' elliptic problems}

\author{Nikolaos Halidias \\
University of the Aegean \\ Department of Statistics and Actuarial
Science
\\ Karlovassi, 83200 \\ Samos \\ Greece \\ email: nick@aegean.gr}

\maketitle


\begin{abstract} In this paper
we are going to show the existence of a nontrivial solution to the
following model problem,
\begin{equation*}
\left\{
\begin{array}{lll}
-\Delta (u)  = 2uln(1+u^2)+\frac{|u|^2}{1+u^2}2u+u(sin(u)-cos(u))
\mbox{ a.e. on } \Omega
\\ \frac{\partial u}{\partial \eta} = 0
 \mbox{ a.e. on } \partial \Omega.
\end{array}
\right.
\end{equation*}
As one can see the right hand side is superlinear. But we can not
use an Ambrosetti-Rabinowitz condition in order to obtain that the
corresponding energy functional satisfies (PS) condition. However,
it follows that the energy functional satisfies the Cerami (PS)
condition.\footnote{2000 Mathematics Subject Calssification:
35A15, 35J20, 35J25
 }
\end{abstract}

\section{Introduction}In this paper
we are going to show the existence of a nontrivial solution to the
following model problem,
\begin{equation}
\left\{
\begin{array}{lll}
-\Delta (u)  = 2uln(1+u^2)+\frac{|u|^2}{1+u^2}2u+u(sin(u)-cos(u))
\mbox{ a.e. on } \Omega
\\ \frac{\partial u}{\partial \eta} = 0
 \mbox{ a.e. on } \partial \Omega.
\end{array}
\right.
\end{equation}
As one can see the right hand side is superlinear. But we can not
use an Ambrosetti-Rabinowitz condition in order to obtain that the
corresponding energy functional satisfies (PS) condition. Let us
recall the well-known Ambrosetti-Rabinowitz condition:

 There
exists some $\theta > 2$ such that
\begin{eqnarray*}
0 < \theta F(u) \leq f(u)u,
\end{eqnarray*}
for all $|u| > M$ for big enough $M$. Here, by $F$ we denote the
function $F(u) = \int_0^r f(r)dr$.

We can see that for our model problem there is not such a $\theta
> 2$.

For such kind of problems there are some papers that extends the
well-known Ambrosetti-Rabinowitz condition. For example one can
also see the very interesting result of D.G. de Figueiredo-J. Yang
\cite{Fig} who considers semilinear problems such that the
corresponding energy functional does not satisfy a (PS) condition.
Also, Gongbao Li-HuanSong Zhou \cite{LIZHOU} made some progress in
this direction. But they assume that $f(u) \geq 0$ for all $u \in
\rr$. Take problem (1) and see that  that $f (u) \to -\infty$ as
$u \to -\infty$, thus we can not say that $f (u) \geq 0$ for all
$u \in \rr$. So, we can not use their method in order to obtain a
nontrivial solution. Finally, let us mention the work of
Costa-Magalhaes \cite{CMAG}. In this paper we extend the results
of \cite{CMAG} using deferent arguments in our proof. The authors
there have proposed the following hypotheses, among others,
\begin{eqnarray*}
f(s)s-p F(s) \geq a |s|^{\mu}, \mbox{ for all } s \in \rr,
\end{eqnarray*}
with $\mu \geq \frac{N/p}{q-p}$ for some $q < p^* =
\frac{Np}{N-p}$. Thus $f(s)s-pF(s)$ must grows faster than
$s^\frac{N/p}{q-p}$.

Our existence theorem, considers more general Neumann problems
than our model problem and at the end of the paper we give a
second example. We must also note that non of the above papers
considers Neumann problems.

Let us introduce the $(PS)$ that we are going to use. \\
{\bf Cerami $(PS)$ condition} Let $X$ be a Banach space and $I:X
\to \rr$.
 For every $\{ u_n \} \subseteq X$ with $|I(u_n)|
 \leq M$ and $(1+\nn u_n \nn)<I^{'}(u_n),\phi> \to 0$ for
 every $\phi \in X$ there exists a strongly
 convergent subsequent. This condition has introduced by Cerami
 (see \cite{Cerami}, \cite{BBF}).

\section{Basic Results}
We are going to show an existence result for the following Neumann
problem,
\begin{equation}
\left\{
\begin{array}{lll}
-\Delta_p (u)  = f(x,u) \mbox{ a.e. on } \Omega
\\ \frac{\partial u}{\partial \eta_p} = 0
 \mbox{ a.e. on } \partial \Omega, \; p \geq 2.
\end{array}
\right.
\end{equation}
We suppose that $\Omega$ is a bounded domain with sufficient
smooth boundary $\partial \Omega$. By $\Delta_p$ we denote the
well-known $p$-Laplacian operator, i.e. $\Delta_p (u) = div(\nn Du
\nn^{p-2}Du)$.

From now on we will denote by $F(x,u) = \int_o^r f(x,r)dr$ and
$h(x,u) = \frac{F(x,u)}{|u|^p}$. We suppose  the following
assumptions on $f$,

{\bf H(f)}$f: \Omega \times \rr \to \rr$ is a Carath\'eodory
function such that
\begin{enumerate}
\item[(i)] for almost all $x \in \Omega$ and for all $u \in \rr$ we have that $|f(x,u)| \leq
C(1+|u|^{\tau-1})$, with $\tau < p^* = \frac{np}{n-p}$ and $C>0$;

\item[(ii)]  there exists some $r < p$ such that for every $k \in
\rr^{+}$ we can find big enough $M>0$ such that
\begin{eqnarray*}
0<(\frac{k}{|u|^r}+p)F(x,u) \leq f(x,u)u
\end{eqnarray*}
 for almost all $x \in \Omega$  for all $u \in \rr$ with
 $|u|>M$;
\item[(iii)] uniformly for almost all $ x \in \Omega$ we have
 $\limsup_{u \to 0} h(x,u) \leq \theta (x)$, with $\theta (x) \leq
 0$ and $\int_{\Omega} \theta (x))dx
 < 0$.
\end{enumerate}

\begin{remark} Note that condition $H(f)(ii)$ is weaker than the
classical condition of Ambrosetti-Rabinowitz. It is easy to see
from $H(f)(ii)$ that $$F(x,u) \geq c
|u|^p\frac{1}{e^{\frac{1}{r|u|^r}}},$$ with $c>0$. From this we
conclude that $\liminf_{|r| \to \infty} h(x,r) > 0$.
\end{remark}

Let us define first the energy functional $I:W^{1,p}(\Omega) \to
\rr$ by $I(u) = \frac{1}{p} \nn Du \nn_p^p - \int_{\Omega}
F(x,u(x))dx$. Under  conditions $H(f)$ it is well known that $I$
is well defined and is a $C^1$ functional. We are going to use the
Mountain-Pass Theorem, so our first lemma
 is that $I$ satisfies the Cerami $(PS)$ condition.

\begin{lemma}
$I$ satisfies the Cerami $(PS)$ condition.
\end{lemma}

{\bf Proof}

Let $\{ u_n \} \subseteq W^{1,p}(\Omega)\cap L^{\infty}(\Omega)$
such that $|I(u_n)| \leq M$ and $(1+\nn u_n
\nn_{1,p})<I^{'}(u_n),\phi> \to 0$ for every $\phi \in
W^{1,p}(\Omega)\cap L^{\infty}(\Omega)$. We must show that $u_n$
is bounded. Suppose that $\nn u_n \nn_{1,p} \to \infty$. We will
show that $\nn D u_n \nn_p \to \infty$. Indeed, from the choice of
the sequence we have

\begin{eqnarray*}
-M \leq \nn D u_n \nn_p^p - \int_{\Omega} p F(x,u_n)dx
\leq M, \Rightarrow \\
\int_{\Omega} F(x,u_n)dx \leq M + \nn D u_n \nn_p^p \Rightarrow \\
c \nn u_n \nn_p^p \leq M + \nn D u_n \nn_p^p,
\end{eqnarray*}
here have used $H(f)(ii)$ (see also Remark 1).

 Thus, we can not suppose that $\nn D u_n
\nn_p$ is bounded because then $\nn u_n \nn_p \to \infty$ and then
from the above relation we obtain a contradiction. So, it follows
that $\nn D u_n \nn_p \to \infty$ and moreover
\begin{eqnarray*}
\nn u_n \nn_p^p \leq (\vep_n + c)\nn D u_n \nn_p^p,
\end{eqnarray*}
with $\vep_n \to 0$.

It follows then that  there exists some $c_1,c_2$ such that

\begin{eqnarray}
c_1 \nn u_n \nn_{1,p} \leq \nn D u_n \nn_p \leq c_2 \nn u_n
\nn_{1,p}.
\end{eqnarray}

Then, from the choice of the sequence it follows
\begin{eqnarray}
-M \leq -\nn Du_n \nn^p_p + \int_{\Omega} p F(x,u_n) dx \leq M,
\end{eqnarray}
and choosing $\phi = u_n$
\begin{eqnarray}
-\vep_n \frac{\nn u_n \nn_{1,p}}{1+\nn u_n \nn_{1,p}} \leq \nn
Du_n \nn_p^p - \int_{\Omega} f(x,u_n) u_ndx \leq \vep_n \frac{\nn
u_n \nn_{1,p}}{1+\nn u_n \nn_{1,p}}.
\end{eqnarray}
Consider now the sequence  $a_n = \frac{1}{ \nn u_n \nn_{1,p}^r}$.
Then multiply inequality (4) with $a_n+1$, substituting with (5)
and using (3) we arrive at
\begin{eqnarray}
& & c \nn u_n \nn_{1,p}^{p-r} \leq  a_n \nn Du_n \nn_p^p \leq \nonumber \\
& & \int_{\Omega}(a_n+1)p F(x,u_n) - f(x,u_n)u_ndx + (a_n+1)M
+\vep_n \frac{\nn u_n \nn_{1,p}}{1+\nn u_n \nn_{1,p}}.
\end{eqnarray}

Let $y_n(x) = \frac{u_n(x)}{\nn u_n \nn_{1,p}}$. Then, it is clear
that there exists some $k \in \rr$ such that $|y_n(x)| \leq k$
a.e. on $\Omega$.

Let $\Omega_1 = \{ x \in \Omega:|u(x)| \leq M \}$. In view of
$H(f)(i)$ we have that for every $M>0$ there exists some $C>0$
such that
\begin{eqnarray*}
\int_{\Omega_1} (a_n+1)p F(x,u_n(x)) - f(x,u_n(x))u_n(x)dx\leq C.
\end{eqnarray*}

Also we have
\begin{eqnarray*}
\int_{\Omega \setminus\Omega_1}a_nF(x,u_n(x))dx+\int_{\Omega
\setminus\Omega_1}p F(x,u_n(x))-f(x,u_n(x))u_n(x)dx \leq
\\
\int_{\Omega \setminus\Omega_1}(\frac{k^r}{|u_n(x)|^r} +p)
F(x,u_n(x)) -f(x,u_n(x))u_n(x)dx.
\end{eqnarray*}

Choose now big enough $M>0$ such that $H(f)(ii)$ holds for $k^r$.
Going back to (6), we obtain a contradiction to the hypothesis
that $u_n$ is not bounded.

In order to show that $I$ satisfies the $(PS)$ condition in
$W^{1,p}(\Omega)$ we have to show that for every $\{ u_n \}
\subseteq W^{1,p}(\Omega)$ with $|I(u_n)| \leq M$ and $(1+\nn u_n
\nn_{1,p})<I^{'}(u_n),\phi> \to 0$ for every $\phi \in
W^{1,p}(\Omega)$ then $\{ u_n \}$ must be bounded in
$W^{1,p}(\Omega)$. But it is well-known that $W^{1,p}(\Omega) \cap
L^{\infty}(\Omega)$ is dense to $W^{1,p}(\Omega)$. So, for every
$u_n$ there exists some sequence $\hat{u}_n^k \subseteq
W^{1,p}(\Omega) \cap L^{\infty}(\Omega)$ such that $\lim_{k \to
\infty} \hat{u}_n^k \to u_n$ strongly in $W^{1,p}(\Omega)$.

From the continuity of $I,I^{'}$ we have that $\hat{u}_n^k$ must
also  satisfy $|I(\hat{u}^k_n)| \leq M$ and $(1+\nn \hat{u}^k_n
\nn_{1,p})<I^{'}(\hat{u}^k_n),\phi> \to 0$ for every $\phi \in
W^{1,p}(\Omega)\cap L^{\infty}(\Omega)$ and for big enough $k \in
\n$.

So, if we suppose that $\nn u_n \nn_{1,p} \to \infty$ then it
follows that $\nn \hat{u}_n^k \nn_{1,p} \to \infty$ and that is a
contradiction following the previous arguments.

Finally, using well-known arguments we can prove that in fact $\{
u_n \}$ have a convergent subsequence.

\QED

\begin{lemma} There exists some $\xi \in \rr$  such
that $I(\xi ) \leq 0$.
\end{lemma}

{\bf Proof}

We claim that there exists big enough $\xi \in \rr$ such that
$I(\xi) \leq 0$.  Suppose not. Then there exists a sequence $\xi_n
\to \infty$ such that $I(\xi_n ) \geq c > 0$. That means
\begin{eqnarray*}
- \int_{\Omega} F(x,\xi_n  )dx \geq c > 0.
\end{eqnarray*}
Using now $H(f)(ii)$ (see also Remark 1) we can say that for
almost all $x \in \Omega$ and all $u \in \rr$ we have that $F(x,u)
\geq \mu |u|^p -c$. So, it follows that
\begin{eqnarray*}
\mu |\xi_n|^p  \leq c \mbox{ for every } n \in \n.
\end{eqnarray*}
But this is a contradiction.

\QED

\begin{lemma} There exists some $\rho >0$ small enough and $a >0$ such that
$I(u) \geq a$ for all $\nn u \nn_{1,p} = \rho$ with $u \in
W^{1,p}(\Omega)$.
\end{lemma}

{\bf Proof} Suppose that this is not true. Then there exists a
sequence $\{ u_n \} \subseteq W^{1,p}(\Omega)$ such as $\nn u_n
\nn_{1,p} = \rho_n$ with $\rho_n \to 0$, with the property that
$I(u_n) \leq 0$. So we arrive at
\begin{eqnarray}
\nn Du_n \nn_p^p \leq p\int_{\Omega} F(x,u_n(x))dx
\end{eqnarray}
 Let $y_n(x) = \frac{u_n(x)}{\nn u_n \nn_{1,p}}$. Using $H(f)(i),(iii)$ we
can prove that there exists $\gamma > 0$ such that
\begin{eqnarray*}
p F(x,u) \leq (\theta (x) + \vep)|u|^p + \gamma |u|^{p^*}
\end{eqnarray*}
Take in account the last estimation and dividing (7) with $\nn u_n
\nn_{1,p}^p$ we arrive at
\begin{eqnarray}
\nn Dy_n \nn_p^p  \leq \int_{\Omega} (\theta (x) +
\vep)|y_n(x)|^pdx + \gamma_1 \nn u_n \nn_p^{p^*-p}.
\end{eqnarray}

Therefore, we obtain that $\nn Dy_n \nn_p \to 0$. Recall that $y_n
\to y$ strongly in $L^p(\Omega)$. Using the lower semicontinuity
of the norm we arrive at $\nn Dy \nn_p \leq \liminf \nn Dy_n \nn_p
\leq \limsup \nn Dy_n \nn_p \to 0$. Thus, $\nn Dy \nn_p = 0$ and
from this we deduce that $y = \xi \in \rr$. But we have that $y_n
\to y$ weakly in $W^{1,p}(\Omega)$ and that $\nn Dy_n \nn_p \to
\nn Dy\nn_p$, so from the uniform convexity of $W^{1,p}(\Omega)$
we obtain that $y_n \to y$ strongly in $W^{1,p}(\Omega)$ and using
the fact that $\nn y_n \nn_{1,p} =1$ we conclude that $\xi \neq
0$.

Going back to (8) and taking the limit we arrive at
\begin{eqnarray*}
\int_{\Omega}  \theta (x)|\xi|^pdx \geq 0.
\end{eqnarray*}
But this is a contradiction.

 \QED

Using now the well-known Mountain Pass theorem we obtain the
desired result.

\section{Applications to Differential Equations}

Consider the following elliptic equation,
\begin{equation}
\left\{
\begin{array}{lll}
-\Delta (u)  = 2uln(1+u^2)+\frac{|u|^2}{1+u^2}2u+u(sin(u)-cos(u))
\mbox{ a.e. on } \Omega
\\ \frac{\partial u}{\partial \eta} = 0
 \mbox{ a.e. on } \partial \Omega.
\end{array}
\right.
\end{equation}
Here, as before, $\Omega \subseteq \rr^n$ is a bounded domain with
smooth enough boundary $\partial \Omega$.

We can check that the corresponding energy functional does not
satisfy an Ambrosetti-Rabinowitz type condition. Moreover, we can
not say that $f(u) \geq 0$ nor that $f(u)+f(-u) = 0$, thus we can
not use the arguments of \cite{LIZHOU}, \cite{LIZHOU2}, even if
our problem had Dirichlet boundary conditions. Finally, $f(\cdot)$
does not satisfy the condition of \cite{Fig} because does not
exist $p>1$ such that $f(u) \geq \mu u^p$ for all $u \geq T$ for
big enough $T$. Also, we can choose big enough $n \in \n$ (i.e.
the dimension of our problem) and see that the above problem does
not satisfy the conditions of \cite{CMAG}.

However, we can check that $f(\cdot)$ satisfies the conditions
that we have proposed.

We can see also that $h$ did not have to go to infinity. Take for
example as $F(u) = |u|^p (a+(b-a)\frac{|u|^r}{1+|u|^r})$. Choose
$a < \frac{\la_1}{p}, b > 0$ and for a suitable choice of $r$ (for
example $r < p$) we can see that $f$ satisfies all the hypotheses
that we have proposed.

\end{document}